\documentclass[12pt,a4paper]{amsart}
\usepackage{amssymb}
\usepackage{multirow}
\usepackage{graphicx}
\usepackage{amsmath,amsbsy,amsfonts,amsthm,latexsym,amsopn,amstext,amsxtra,euscript, amscd,mathrsfs}
\usepackage[noadjust]{cite}

\theoremstyle{plain}
\newtheorem{theorem}{Theorem}[section]

\newtheorem{conjecture}{Conjecture}

\theoremstyle{definition}

\theoremstyle{remark}

\numberwithin{equation}{section}

\numberwithin{table}{section}

\numberwithin{figure}{section}

\begin{document}

\title[An observation on the determinant of a Sylvester-Kac type matrix]{An observation on the determinant of a Sylvester-Kac type matrix}

\author{Carlos M. da Fonseca}
\address{Kuwait College of Science and Technology, Doha District, Block 4,
		P.O. Box 27235, Safat 13133, Kuwait}
\email{c.dafonseca@kcst.edu.kw}
\address{University of Primorska, FAMNIT, Glagoljsa\v ska 8, 6000 Koper, Slovenia}
\email{carlos.dafonseca@famnit.upr.si}
\author{Emrah K\i l\i \c{c}}
\address{TOBB University of Economics and Technology, Mathematics
		Department, 06560 Ankara, Turkey}
\email{ekilic@etu.edu.tr}
\date{\today }
\begin{abstract}
Based on a less-known result, we prove a recent conjecture concerning the determinant of a certain Sylvester-Kac type matrix and consider an extension of it. 
\end{abstract}

\maketitle

\section{The Conjecture}

Quite recently, in order to find formulas for the determinants of some Lie algebras,  Z. Hu and P.B. Zhang proposed in \cite{HZ2019} the following conjecture.
\begin{conjecture} \label{conj}
The determinant of the matrix $(n+1)\times (n+1)$ 
\begin{equation*}\label{SK1}
J_{n}(z_0,z_1)=\left(
	\begin{array}{cccccc}
	z_0+nz_1 & 1 &  &  &  &  \\
	n & z_0+(n-2)z_1 & 2 &  &  &  \\
	& n-1 & z_0+(n-4)z_1 & \ddots  &  &  \\
	&  & \ddots  & \ddots  & n-1 &  \\
	&  &  & 2 & z_0-(n-2)z_1 & n \\
	&  &  &  & 1 & z_0-nz_1 
	\end{array}
	\right) \,
\end{equation*}
is 
$$
\prod_{k=0}^{n}\left( z_0-(n-2k)\sqrt{z_1^2+1}\right)\, .
$$
\end{conjecture} 

Notice that Conjecture \ref{conj} is equivalent to state that the eigenvalues of  $J_{n}(0,z_1)$ are 
$$
\pm (n-2k)\sqrt{z_1^2+1}\, , \quad \mbox{for $k=0,1,\ldots,\lfloor n/2\rfloor$.}
$$

The matrix $J_{n}(z_0,z_1)$ can be easily identified as an extension of the so-called Sylvester-Kac matrix. In fact, setting $z_1=0$ we find the characteristic matrix of the  Sylvester-Kac matrix, also known as Clement matrix, 
$$\left(
\begin{array}{ccccc}
0 & 1 &&& \\
n & 0 & 2 && \\
& n-1 & 0 & \ddots & \\
&& \ddots & \ddots & n \\
&&& 1 & 0 \\
\end{array}
\right)\, .
$$ 
The characteristic polynomial of this  matrix (that is, $\det J_{n}(x,0)$) was first conjecture in \cite{S1854}, by the $19$th century British mathematician James Joseph Sylvester celebrated,  among other facets, as the founder of the American Journal of Mathematics, in $1878$.

A fully comprehensive list of results on the different proofs for Sylvester's conjecture and the eigenpairs of non-trivial extensions of the Sylvester-Kac matrix  can be found in \cite{A2005,BR2007,C2010,CW2008,C1959,EK1994,dFMMW2013,FS1965, H2005,I2002, K1947,K2013,KA2016,M1923,OJ2017,R1957,S1926,TT1991,V1964}. 

The aim of this short note is to prove Conjecture \ref{conj} based on a result by W. Chu in \cite{C2010}. We also provide a general result containing other particular known determinants.

\section{An extension to Sylvester-Kac matrix}

In $2010$, cleverly based  on two generalized Fibonacci sequences, W. Chu proved the following theorem.

\begin{theorem}[\cite{C2010}] \label{TheoremGeneral}
	The determinant of the matrix $(n+1)\times (n+1)$ 
	\begin{equation*}\label{SK1}
	M_{n}(x,y,u,v)=\left(
	\begin{array}{cccccc}
	x & u &  &  &  &  \\
	nv & x+y & 2u &  &  &  \\
	& (n-1)v & x+2y & \ddots  &  &  \\
	&  & \ddots  & \ddots  & n-1 &  \\
	&  &  & 2v & x-(n-1)y & nu \\
	&  &  &  & v & x+ny 
	\end{array}
	\right) \,
	\end{equation*}
	is 
	$$
	\prod_{k=0}^{n}\left( x+\frac{ny}{2}+\frac{n-2k}{2}\sqrt{y^2+4uv}\right)\, .
	$$
\end{theorem}  

Of course, the formula for the determinant in Theorem \ref{TheoremGeneral} can be rewritten as
$$
\prod_{k=0}^{\lfloor n/2\rfloor}\left(\left( x+\frac{ny}{2}\right)^2-\frac{(n-2k)^2}{4}(y^2+4uv)\right)\, .
$$

Now setting $x=z_0+nz_1$, $y=-2z_1$, and $u=v=1$, we prove immediately Conjecture \ref{conj}. 

Moreover, in the spirit of \cite{A2005,dFMMW2013,H2005}, using Theorem \ref{TheoremGeneral}, we can also conclude the following theorem.

\begin{theorem}
The eigenvalues of  
	\begin{equation*}
	M^\pm_n(a,b,r) = \left( \begin{array}{ccccccc}
	n a r & b &&&&&\\
	na & ((n-1)a \pm b)r & 2b &&&& \\
	& (n-1)a & ((n-2)a\pm 2b)r &  3b &&& \\
	&& (n-2)a & \ddots & \ddots  && \\
	&&& \ddots & \ddots & nb &\\
	&&&& a & \pm nb r & \\
	\end{array} \right) \label{M} .
	\end{equation*}
are
$$\frac{1}{2} \left(n r(a\pm b)+(n-2k)\sqrt{4ab+r^2(a\mp b)^2}\right)\, ,$$ for $k=0,1,\ldots,n$....
\end{theorem}

\end{document}